\documentclass{amsart}
\usepackage{amsmath}
\usepackage{amssymb}

\begin{document}

\author[S.Albeverio, V.Koshmanenko, M.Pratsiovytyi, G.Torbin]
{Sergio Albeverio$^{1,2,3,4}$, \thinspace Volodymyr
Koshmanenko$^5$,
 Mykola Pratsiovytyi$^6$, Grygoriy Torbin$^7$}
\title[]{$\widetilde{Q}-$ representation of real numbers and fractal probability
distributions }
\maketitle

\begin{abstract}
A $\widetilde{Q}-$representation of real numbers is introduced as
a generalization of the $s-$adic expansion. It is shown that the
$\widetilde{Q}-$representation is a convenient tool for the
construction and study of fractals and measures with complicated
local
structure. Distributions of random variables $\xi $ with independent $%
\widetilde{Q}-$symbols are studied in details. Necessary and
sufficient conditions for the corresponding probability measures
$\mu _{\xi }$  to be either absolutely continuous or singular
(resp. pure continuous, or pure point) are found in terms of the
$\widetilde{Q}-$representation. The metric, topological, and
fractal properties for the distribution of $\xi $ are
investigated. A number of examples are presented.
\end{abstract}

$^1$I{nstitut f\"{u}r Angewandte Mathematik, Universit\"{a}t Bonn,
Wegelerstr. 6, D-53115 Bonn\ (Germany); }$^2${SFB 611, \ Bonn, \
BiBoS, (Bielefeld - Bonn); }$^3${IZKS Bonn}; $^4${CERFIM, Locarno
and Acc. Arch. (USI) (Switzerland).} {E-mail: albeverio@uni-bonn.de}

$^5$Institute of Mathematics, Tereshchenkivs'ka str. 3, 01601 Kyiv
(Ukraine) , \ {E-mail: kosh@imath.kiev.ua}

$^6$National Pedagogical University, Pyrogova str. 9, 01030 Kyiv
(Ukraine), \ {E-mail: prats@ukrpost.net}

$^7$National Pedagogical University, Pyrogova str. 9, 01030 Kyiv
(Ukraine), {E-mail: torbin@imath.kiev.ua}

\medskip \textbf{AMS Subject Classifications (2000):11A67, 11K55, \thinspace
26A30, 26A46, 28A80, \thinspace 28D05, 60G30}

\textbf{Key words: } $\widetilde{Q}$-representation of real
numbers; random variables with independent
$\widetilde{Q}$-symbols; Hausdorff dimension, fractals; pure
point, singular continuous, absolutely continuous measures; C, P,
S types of measures.

\section{Introduction}

As well known there exist only three types of pure probability
distributions: discrete, absolutely continuous and singular.
During a long period mathematicians had a rather low interest in
singular probability distributions, which was mainly caused by the
two following reasons: the absence of effective analytic tools and
the widely spread point of view that such distributions do not
have any applications, in particular in physics, and are
interesting only for theoretical reasons. The interest in singular
probability distributions increased however in 1990's due their
deep connections with the theory of fractals. On the other hand,
recent investigations show that singularity is generic for many
classes of random variables, and absolutely continuous and
discrete distributions arise only in exceptional cases (see, e.g.
\cite{P,Z}). Possible applications in the spectral theory of
self-adjoint operators (\cite{Tr}) is an additional reason in the
intensive investigation of singular continuous measures. It was
proven that Schr\"odinger type operators with singular continuous
spectra  are generic for some classes of potentials (\cite{DJMS}).
Moreover, by using the fractal analysis of the corresponding
spectral singular continuous measures, it is possible to analyze
the dynamical properties of the corresponding quantum systems
(\cite{L}).

Usually the singular probability distributions are associated with
the Cantor-like distributions. Such distributions are supported by
nowhere dense sets of zero Lebesgue measures. In the sequel we shall
call such distributions the distributions of the {\it C-type}. But
there exist singular probability distributions with other
metric-topological properties of their topological support.

A set $E \subseteq R^1$ is said to be  of the pure {\it S-type} if
it is a union of an at most countable family of closed intervals $[
a_i, b_i] , \ \ a_i < b_i$. ~~~~~A set $E \subseteq R^1 $ is said to
be of the pure {\it C-type} if it is a nowhere dense set of zero
Lebesgue measure. ~~~~ A set $E \subseteq R^1 $ is said to be of the
pure {\it P-type} if it is a nowhere dense set such that for any $x
\in E$ and for any $\varepsilon
>0$ the set $E \bigcap (x-\varepsilon, x+\varepsilon)$ is of
positive Lebesgue measure. ~~~~~

A probability measure $ \mu $ on $R^1$ is said to be of the pure
S-type if the measure is supported by a  set of the pure S-type,
i.e., if there exists a set $S$ of the pure S-type such that the
topological support (the minimal closed set supporting the measure)
$S _{\mu }=S^{cl}$ and $ \mu(S)=1$, where $cl$ stands for the
closure.

 A probability measure $ \mu $ on $R^1$ is said to be of the pure C-type
 if its topological support $S _{\mu }$ is a nowhere dense set
 and the measure is supported by a set  of the pure C-type, i.e.,  there
exists a  set $C$ of C-type such that $S _{\mu }=C^{cl}$ is a
nowhere dense set and $ \mu(C)=1$.

A probability measure $ \mu $ on $R^1$ is said to be of the pure
P-type if its topological support $S _{\mu }$ is a nowhere dense set
and the measure is supported by a set  of the pure P-type, i.e., if
there exists a set $P$ of the pure  P-type such that $S _{\mu
}=P^{cl}$ is a nowhere dense set and $\mu(P)=1$.

In \cite{P} it has been proved that any \textit{singular continuous}
probability measure $\mu_{\rm sing}$ on $R^1$ can be represented in
the following form:
\begin{equation*}
\mu_{\rm sing}= \alpha_1\mu_{\rm sing}^C+ \alpha_2\mu_{\rm sing}^S+
\alpha_3\mu_{\rm sing}^P,
\end{equation*}
where $\alpha_i\geq 0, \alpha_1 + \alpha_2 + \alpha_3
=1,~~~~\mu_{{\rm sing}}^C$, $\mu_{{\rm sing}}^S$  and $\mu_{{\rm
sing}}^P$ are singular continuous probability measures of the pure
$C-$, $S-$ resp. $P-$types.
 A corresponding decomposition holds for the singular component of
 the spectral measure  and the Hilbert space for any given self-adjoint operator with a simple spectrum
 (see \cite{AKT} for details).

In the paper we introduce into consideration the so-called
$\widetilde{Q}$-representation of real numbers which is a
convenient tool for the construction of a wide class of fractals.
This class contains  Cantor-like sets as well as everywhere dense
noncompact fractals with any desirable Hausdorff-Besicovitch
dimension $\alpha_0 \in [0,1]$. By using the $\widetilde{Q}$
-representation we introduce a family of random variables with independent $%
\widetilde{Q}$-symbols. This family contains all possible the
above mentioned types of singular continuous measures, and (as a
very particular case) the class of all self-similar measures on
$[0,1]$ satisfying the open set condition (see section 4 for
details).

An additional reason for the investigation of the distribution of
the random variables with independent $\widetilde{Q}$-symbols is
to extend the so-called Jessen-Wintner theorem to the case of sums
of random variables which are not independent. In fact this
theorem asserts that if a random variable is a sum of the
convergent series of  independent discretely distributed random
variables, then it has a pure distribution. In general setting,
necessary and sufficient conditions for such probability
distributions to be singular resp. absolutely continuous are still
unknown.

In this paper we completely investigated the structure of the
random variables with independent $\widetilde{Q}$-symbols. The
topological, metric and fractal properties of the above mentioned
 probability distributions are also studied.

\section{\protect\bigskip $\widetilde{Q}$-representation of real numbers}

We describe the notion of the $\widetilde{Q}-$\textit{%
representation }for real numbers $x\in \left[ 0,1\right] $. Let us consider a $\mathbf{N}_k\times \mathbf{N}-$%
matrix $ \widetilde{Q}=\left\| q_{ik}\right\|,$ $ i \in
\mathbf{N}_k$, $k\in \mathbf{N}$, where $\mathbf{N}$ stands for the set of natural numbers and $\mathbf{N}%
_k=\{0,1,...,{N}_k\}$, with $0<{N}_k\leq \infty $.  We suppose
that
\begin{equation}
q_{ik}>0 ~~~ , ~ i \in \mathbf{N}_k, ~ k\in \mathbf{N}. \label{I}
\end{equation}
Besides, we assume that for each $k\in \mathbf{N}$:
\begin{equation}
\sum_{i\in \mathbf{N}_k}q_{ik}=1,  \label{1}
\end{equation}
and
\begin{equation}
\prod_{k=1}^\infty \max_{i\in \mathbf{N}_k}\left\{ q_{ik}\right\} =0.
\label{2}
\end{equation}

Given a $\widetilde{Q}-$matrix we consecutively perform decompositions of
the segment $\left[ 0,1 \right]$ as follows.

Step 1. We decompose $\left[ 0,1 \right] $ (from the left to the
right) into the union of closed intervals
$\Delta^{\widetilde{Q}}_{i_1}$, $i_1\in \mathbf{N}_1$ (without
common interior points) of the length $\left|
\Delta^{\widetilde{Q}}_{i_1}\right| =q_{i_11}$,
\begin{equation*}
\left[ 0,1 \right] =\bigcup_{i_1 \in \mathbf{N}_1}
\Delta^{\widetilde{Q}}_{i_1}.
\end{equation*}

Each interval $\Delta^{\widetilde{Q}}_{i_1}$ is called a {\it 1-rank
interval}.



Step $k\geq2$. We decompose (from the left to the right) each closed
$(k-1) $-rank interval $\Delta^{\widetilde{Q}}_{i_1i_2...i_{k-1}}$
into the union of closed {\it $k-$rank intervals}
$\Delta^{\widetilde{Q}}_{i_1i_2...i_k},$
\begin{equation*}
\Delta^{\widetilde{Q}}_{i_1i_2...i_{k-1}}=\bigcup_{i_k \in \mathbf{N}_k}%
\Delta^{\widetilde{Q}}_{i_1i_2...i_k},
\end{equation*}
where their lengths
\begin{equation}
\left| \Delta^{\widetilde{Q}}_{i_1i_2...i_k}\right|= q_{i_11}\cdot
q_{i_22}\cdot \cdot \cdot q_{i_kk}=\prod_{s=1}^kq_{i_ss}  \label{3}
\end{equation}
are related as follows
\begin{equation*}
\left| \Delta^{\widetilde{Q}}_{i_1i_2...i_{k-1}0}\right| :\left|
\Delta^{\widetilde{Q}}_{i_1i_2...i_{k-1}1}\right| :\cdot \cdot \cdot
:\left| \Delta^{\widetilde{Q}}_{i_1i_2...i_{k-1} i_k}\right| :\cdot
\cdot \cdot =q_{0k}:q_{1k}:\cdot \cdot \cdot :q_{i_k k}:\cdot \cdot
\cdot~~~~~~~~.
\end{equation*}

For any sequence of indices $\{ i_k \},$  $i_k\in \mathbf{N}_k$,
there corresponds the sequence of embedded closed intervals
\begin{equation*}
\Delta^{\widetilde{Q}}_{i_1}\supset \text{ }\Delta^{\widetilde{Q}}
_{_{i_1i_2}}\supset \cdot \cdot \cdot \supset \Delta^{\widetilde{Q}}
_{_{i_1i_2...i_k}}\supset \cdot \cdot \cdot
\end{equation*}
such that $|\Delta^{\widetilde{Q}}_{i_1...i_k}|\to 0$, $k\rightarrow
\infty $, due to (\ref {2}) and (\ref{3}). Therefore, there exists a
unique point $x\in \left[ 0,1\right] $ belonging to all intervals
$\Delta^{\widetilde{Q}}_{i_1},$ $\Delta^{\widetilde{Q}}_{i_1i_2},$
..., $\Delta^{\widetilde{Q}}_{i_1i_2...i_k},...$. Conversely, for
any point $x\in \left[ 0,1\right] $ there exists a sequence of
embedded intervals $\Delta^{\widetilde{Q}}_{i_1}\supset
\Delta^{\widetilde{Q}}_{i_1i_2}\supset ...\supset $
$\Delta^{\widetilde{Q}}_{i_1i_2...i_k}\supset ...$ containing $x$,
i.e.,
\begin{equation}
x=\bigcap_{k=1}^\infty
\Delta^{\widetilde{Q}}_{i_1i_2...i_k}=\bigcap_{k=1}^\infty
\Delta^{\widetilde{Q}}_{i_1(x)i_2(x)...i_k(x)}=:\Delta^{\widetilde{Q}}_{i_1(x)i_2(x)...i_k(x)...}
\label{4}
\end{equation}
Notation (\ref{4}) is called the $%
\widetilde{Q}-$\textit{representation} of the point $x\in \left[
0,1\right] . $

\textbf{Remark 1.} The correspondence $\left[ 0,1 \right] \in x
\Leftrightarrow \{(i_1(x), i_2(x), ..., i_k(x), ... )\} $ in
(\ref{4}) is one-to-one , i.e., the $\widetilde{Q}-$representation
is unique for every point $x \in \left[ 0,1 \right] $, provided
that the $\widetilde{Q}$-matrix contains an infinite number of
columns with an infinite number of elements. However in the case,
where $N_k<\infty,$ $\forall k>k_0$ for some $k_0$, there exists a
countable set of points $x\in \left[ 0,1\right] $ having two different $%
\widetilde{Q}-$representations. Precisely, this is the set of all
end-points of intervals $\Delta^{\widetilde{Q}}_{i_1i_2...i_k}$ with
$k>k_0.$

One has the formula
\begin{equation}
x =D_1(x) +\sum_{k=2}^\infty \left[
D_k(x)\prod_{s=1}^{k-1}q_{i_s(x)s} \right] =\sum_{k=1}^\infty D_k(x)
L_{k-1}(x)  \label{5}
\end{equation}
where $ D_k(x):=\left\{
\begin{array}{lll}
0, & \text{if } & i_k(x)=0, \\ \sum\limits_{i=0}^{i_k(x)-1}q_{ik},
& \text{if } & i_k(x) \geq 1
\end{array}
\right.$ and where we put
\begin{equation*}
L_{k-1}(x):=\mid\Delta^{\widetilde{Q}}_{i_1(x)...i_{k-1}(x)}\mid=\prod_{s=1}^{k-1}q_{i_s(x)s}
\end{equation*}
for $k>1$, and $L_{k-1}(x)=1$, if $k=1$.

We note that (\ref{5}) follows from (\ref{4}) since the common
length of all intervals lying on the left side of a point
$x=\Delta^{\widetilde{Q}}_{i_1(x)...i_k(x)...}$ can be calculate as
the sum of all $1-$rank intervals lying on the left from $x$ (it is
the first term $D_1(x)$ in (\ref{5})), plus the sum of all $2-$rank
intervals from $\Delta^{\widetilde{Q}}_{i_1(x)}$, lying on the left
side from $x$ (the second term $D_2(x)\cdot q_{i_1(x)1} $ in
(\ref{5})~), and so on.

\textbf{Remark 2.} If $q_{ik}=q_i\,\,, k\in \mathbf{N},$ then the
\smallskip $\widetilde{Q}-$representation coincides with the $Q-$%
representation (see \cite{TuP}); moreover, if $q_{ik}=\frac 1s$
for some natural number$\,s>1,$ then the $\smallskip \widetilde{Q}-$%
representation coincides with the classical $s-$adic expansion.

\section{$\widetilde{Q}(\mathbf{V})-$representation for fractals}

The $\widetilde{Q}$-representation  allows to construct in a
convenient way a wide class of fractals on $R^1$ and other
mathematical objects with fractal properties. Firstly we consider
compact fractals from $R^1$. Let $\mathbf{V}:=\{\mathbf{V}_k \}_{k=1}^\infty$, $\mathbf{V}_k \subseteq%
\mathbf{N}_k,$ and let us consider the set
\begin{equation}
\Gamma _{\widetilde{Q}(\mathbf{V})}\equiv\Gamma :=\left\{ x \in[0,1%
]:~~ x=\Delta^{\widetilde{Q}}_{i_1i_2...i_k...}, ~ i_k \in
\mathbf{V}_k\right\}, \label{Gamma}
\end{equation}
i.e., $\Gamma$ consists of points, which can be $\widetilde{Q}-$%
represented by using only symbols $i_k$ from the set $\mathbf{V}_k$ on each $%
k$-th position of their $\widetilde{Q}-$representation.

If $\mathbf{V}_k \neq \mathbf{N}_k$ at least for one $k<k_0$, and $\mathbf{V}%
_k=\mathbf{N}_k$ for all $k \geq k_0$ with some fixed $k_0>1$,
then $ \Gamma $ is a union of closed intervals. In this case one
can get $\Gamma $ removing from $[0,1]$ all
open intervals $\dot \Delta^{\widetilde{Q}}_{i_1...i_k},$ $k<k_0$ with $i_k \notin \mathbf{V}%
_k$ (where a point over $\Delta $ means that an interval is open).
If the condition $\mathbf{V}_k \neq\mathbf{N}_k$ holds for
infinitely many values of $k$, then obviously $\Gamma $ is a nowhere
dense set.

We shall study the metric properties of the sets
$\Gamma_{\widetilde{Q}(\mathbf{V})}$. Let $S_k(\mathbf{V})$ denote
the sum of all elements $q_{ik}$ such that $i_k \in \mathbf{V}_k$,
i.e.,
\begin{equation*}
S_k(\mathbf{V}):=\sum_{i\in \mathbf{V}_k}q_{ik}.
\end{equation*}
We note that $0<S_k(\mathbf{V}) \leq 1$ due to (\ref{I}), (\ref{1}).

\textbf{Lemma 1.} \textit{ Lebesgue measure $\lambda (\Gamma )$ of
the set $\Gamma $  is equal to}
\begin{equation}
\lambda (\Gamma )=\prod_{k=1}^\infty S_k(\mathbf{V}).
\label{modul}
\end{equation}

\proof Let $\Gamma _n:=\bigcup\limits_{i_k\in \mathbf{V}_k}\Delta
_{i_1...i_n}$. It is easy to see that $\Gamma _n\subseteq \Gamma
_{n-1}$ and ~~~~$\Gamma =\bigcap\limits_{n=1}^\infty \Gamma _n$.
From the definition of the sets $\Gamma_n$ and from (\ref{3}), it
follows that $\lambda (\Gamma
_n)=\prod\limits_{k=1}^nS_k(\mathbf{V})$, and, therefore, $\lambda
(\Gamma )=\lim\limits_{n\longrightarrow \infty }\lambda (\Gamma _n)=
\prod\limits_{k=1}^\infty S_k(\mathbf{V}).$

\endproof

\textbf{Lemma 2.} \textit{Let $W_k(\mathbf{V})=1-S_k(\mathbf{V})
\geq 0$. The set $\Gamma $ is of  zero Lebesgue measure if and only
if}
\begin{equation}
\sum^{\infty }_{k=1}W_k(\mathbf{V})=\infty ~~~. \text{ }  \label{WK}
\end{equation}
\proof This assertion is a direct consequence of the previous
lemma and the well known relation between infinite products and
infinite series. Namely, for a sequence $0 \le a_k < 1 ,$ the
product $\prod\limits^{\infty}_{k=1 }(1-a_k)=0$ if and only if the
sum $\sum\limits^{\infty}_{k=1 }a_k=\infty$. In our case $a_k =
1-S_k(\mathbf{V})$.
\endproof

The above mentioned procedure allows to construct nowhere dense
compact fractal sets $E$ with desirable Hausdorff-Besicovitch
dimension (including the anomalously fractal case ($\alpha_0
(E)=0$) and the superfractal case ($\alpha_0 (E)=1$)) in a very
compact way.

\textbf{Theorem 1.} {\it Let $\mathbf{N}_k = N^0_{s-1} := \{
0,1,..., s-1 \}$ $ k \in N$, let $\mathbf{V_0}= \{ v_1, v_2, ...,
v_m\}\subset N^0_{s-1}$ and let the matrix $\widetilde{Q}$ have
the following asymptotic property: $$ \lim_{k \to \infty} q_{ik} =
q_i, \,\,\, i \in N^0_{s-1}.$$ Then:

1) the Hausdorff-Besicovitch dimension of the set
$\Gamma_{\widetilde{Q}( \mathbf{V_0})}$ is the root of the
following equation

\begin{equation}\label{falk}
 \sum_{i \in \mathbf{V_0}} q_i^x =1, \,\,\,\, ;
\end{equation}

2) if $$M[\widetilde{Q}, (\nu_0, ..., \nu_{s-1})] = \left\{ x: \,
\Delta^{\widetilde{Q}}_{\alpha_1 (x)...\alpha_k(x)...}, \,\,
\lim_{k \to \infty} \frac{N_i (x,k)}{k} = \nu_i, \, i \in
N^0_{s-1} \right\},$$ where $N_i (x,k) $ is the number  of symbols
''$i$'' among the first $k$ symbols of the
$\widetilde{Q}$-representation of $x$, then

\begin{equation}\label{hb}
  \alpha_0 (M[\widetilde{Q}, (\nu_0, ..., \nu_{s-1})]) = \frac{ \sum\limits_{i=0}^{s-1} \nu_i \ln
\nu_i}{   \sum\limits_{i=0}^{s-1} \nu_i \ln q_i}.
\end{equation}

}

 \proof Firstly we consider the particular  case where the
matrix $\widetilde{Q}$ has  exactly $s$ rows and all its columns
are the same: $q_{ik}=q_i$. In such a simple case the
$\widetilde{Q}$-representation reduces to the $Q$-representation
studied in \cite{P}. It is easy to prove (see, e.g., \cite{TuP}),
that to calculate the Hausdorff-Besicovitch dimension of any
subset $E \subset [0,1]$ it is sufficient  to consider a class of
cylinder sets of different ranks generated by $Q$-partitions of
the unit interval.
 The
Billingsley theorem (see, e.g., \cite{Bi}, p.141) admits a
generalization to the class of $Q$-cylinders, and, from this
theorem it follows that in the case of usual $Q$-representation,
the Hausdorff-Besicovitch dimension of the set $M[Q, (\nu_0, ...,
\nu_{s-1})]$ is equal to the right-side expression in (\ref{hb}).

In the  $Q$-case the set $\Gamma_{Q(\mathbf{V_0})}$ is a
self-similar set satisfying the open set condition. Therefore, the
Hausdorff-Besicovitch dimension of this set is the root of
equation (\ref{falk}).

Now let us consider a general case of theorem 1. To this end we
introduce into consideration the following transformation $f$ of
$[0,1]:$
$$ f(x)= f (\Delta^{Q}_{\alpha_1 (x) ...\alpha_k (x)...})
= \Delta^{\widetilde{Q}}_{\alpha_1 (x)...\alpha_k (x)...~~~.}$$

It is  not hard to prove, that such a transformation  belongs to the
DP-class (see, e.g., \cite{APT}), i.e., $f$ preserves the
Hausdorff-Besicovitch dimension of any subset of $[0,1]$.  Since $
f(\Gamma_{Q(\mathbf{V_0})})= \Gamma_{\widetilde{Q}(\mathbf{V_0})}$
and $ f(M[Q, (\nu_0, ..., \nu_{s-1})])= M[\widetilde{Q}, (\nu_0,
..., \nu_{s-1})]$, we get the desired formulas under general
assumptions of theorem 1.
\endproof
\textbf{ Example 1}. If  $\mathbf{N}_k=\{ 0,1,2\}
,\mathbf{V}_k=\{0,2\},
q_{1k}\rightarrow 0,$ but $\sum\limits_{k=1}^\infty q_{1k}=\infty$ with $%
q_{0k}=q_{2k}=\frac{1-q_{1k}}2,$ then $\Gamma  $ is a nowhere
dense set of zero Lebesgue measure. From Theorem 1 it follows that
the Hausdorff dimension of this set is equal 1. In the terminology
of \cite{P} a set of this kind is called a superfractal set.

\textbf{ Example 2}. If $\mathbf{N}_k=\{ 0,1,2\} ,\mathbf{V}_k=\{ 0,2\}$, $%
q_{1k}\rightarrow 1$ (but $\prod\limits_{k=1}^{\infty}q_{1k}=0$), and $%
q_{0k}=q_{2k}=\frac{1-q_{1k}}2,$ then $\Gamma  $ is a nowhere
dense set of zero Lebesgue measure and of zero Hausdorff
dimension, i.e., $\Gamma $ is an anomalously fractal set (see
\cite{P}).

\section{Random variables with independent $\widetilde{Q}-$symbols}

Let $\{ \xi_k,~~~  k \in \mathbf{N} \}$ be a sequence of independent
random  variables with the following distributions
\begin{equation*}
P(\xi_k=i):=p_{ik}\geq 0, ~~~with~~~\sum_{i \in
\mathbf{N}_k}p_{ik}=1,  ~~ k \in \mathbf{N}.
\end{equation*}

By using $\xi_k$ and the $\widetilde{Q}$-representation we
construct a random variable $\xi$ as follows:
\begin{equation}
\xi: =\Delta^{\widetilde{Q}}_{\xi _1\xi _2...\xi _k...}~.
\label{ksi}
\end{equation}

The distribution of $\xi $ is completely determined by two matrices: $%
\widetilde{Q}$ and $ \widetilde{P}=||p_{ik}||,$  where some
elements of the matrix $\widetilde{P}$ possibly are equal to zero.
Of course, all sets $\mathbf{N}_k$ are the same as those in the
$\widetilde{Q}$-matrix. Let $\mu_\xi $ be the measure
corresponding the distribution of the random variable $\xi$ with
independent $\widetilde{Q}$-symbols.

If $q_{ik}=q_{i}$ and $p_{ik}=p_{i}$ $\forall j\in N$, $i\in
N_{s-1}^0$ (i.e., $\xi$ is a random variable with independent
identically distributed Q-digits), then the measure $\mu_{\xi}$ is
the self-similar measure associated with the list $(S_1, ...,
S_{s-1}, p_1, ...,p_{s-1})$, where $S_i$ is the similarity with
the ratio $q_i$ ($\sum\limits_{i=0}^{s-1} q_{i}=1$), and the list
$(S_1, ..., S_{s-1})$ satisfies the open set condition. More
precisely, $\mu_{\xi}$ is the unique Borel probability measure on
$[0,1]$ such that

$$\mu_{\xi} = \sum_{i=0}^{s-1} p_{i}\cdot \mu_{\xi}\circ S_i^{-1},$$
 (see, e.g., \cite{Fa} for details).
 In the  so-called "$Q^*-$ case" we construct the measure $\mu_{\xi}$ in a
 similar way but with the possibility of changing of the ratios and probabilities from the list
  $(S_1, ..., S_{s-1}, p_1, ...,p_{s-1})$ at each
  stage of the construction.
     In our general "$\widetilde{Q}-$ case" we may  additionally choose the number
  of contracting similarities (including a countable number) at each
  stage of the construction.

From (\ref{5}) it follows that $ \xi $ can be represented as a sum
of the convergent series of discretely distributed random
variables which are not independent. Nevertheless the distribution
of $\xi$ is of pure type.




\textbf{Theorem 2.} \textit{The measure $\mu_\xi $ is of pure
type, i.e., it is either purely absolutely continuous, resp.,
purely point, resp., purely singular continuous. Precisely, }

\textit{1) $\mu_\xi$ is purely absolutely continuous if and only
if }

\textit{
\begin{equation}
\rho := \prod^{\infty }_{k=1}\left\{ \sum _{i\in \mathbf{N}_k}
\sqrt{p_{ik}\cdot q_{ik} }\right\} >0 ;  \label{1)}
\end{equation}
}

\textit{2) $\mu_\xi$ is purely point if and only if }

\textit{
\begin{equation}
P_{max} := \prod_{k=1}^{\infty } \max_{i\in \mathbf{N}_k}\{ p_{ik}
\} >0;  \label{2)}
\end{equation}
}

\textit{3) $\mu_\xi$ is purely singular continuous if and only if
}

\textit{
\begin{equation}
\rho=0=P_{max}.  \label{3)}
\end{equation}
}

\proof Let $\Omega_k = \mathbf{N}_k $, $\mathcal{A}_k =
2^{\Omega_k}$. We define measures $\mu_k$ and $\nu_k$ in the
following way: $$   \,\,\, \mu_k (i) = p_{ik}; \,\,\, \nu_k (i)=
q_{ik}, ~~i \in \Omega_k.$$ Let $$ (\Omega, \mathcal{A}, \mu)=
\prod_{k=1}^\infty (\Omega_k , \mathcal{A}_k, \mu_k ), \,\,\, \,\,
(\Omega, \mathcal{A}, \nu)= \prod_{k=1}^\infty (\Omega_k,
\mathcal{A}_k, \nu_k)$$ be the infinite products of probability
spaces, and let us consider the measurable mapping $f: \Omega \to
[0;1]$ defined as follows: $$ \forall \omega = ( \omega_1,
\omega_2, ..., \omega_k, ...) \in \Omega, \,\,\,\, f(\omega)= x=
\Delta_{i_1 (x) i_2 (x)... i_k(x)... }$$ with $\omega_k = i_k (x)$
$k \in N$.

We define the measures $\mu^*$ and $\nu^*$ as the image measure of
$\mu$ resp. $\nu$ under $f$: $$ \mu^* (B):= \mu (f^{-1} (B)) ;
\,\,\, \nu^*(B)= \nu (f^{-1}(B)),  B \in \mathcal{B}.$$

It is easy to see that $\nu^*$ coincides with Lebesgue measure
$\lambda$ on $[0,1]$, and $\mu^* \equiv \mu_\xi$. In general, the
mapping $f$ is not bijective, but there exists a countable set
$\Omega_0$ such that $\nu(\Omega_0)=0, \mu(\Omega_0)=0 $ and the
mapping $f: \Omega \setminus \Omega_0 \rightarrow [0,1]$ is
bijective.

Therefore, the measure $\mu_\xi$ is absolutely continuous (singular)
with respect  to Lebesgue measure if and only if the measure $\mu$
is absolutely continuous (singular) with respect to the measure
$\nu$. Since, $q_{ik}>0$, we conclude that $\mu_k \ll \nu_k$,
$\forall k \in N$. By using Kakutani's theorem \cite{KAK}, we have
\begin{equation}
 \mu_\xi \ll \lambda \,\,\, \Leftrightarrow \,\,\,
\prod_{k=1}^\infty \int_{\Omega_k} \sqrt{\frac{ d \mu_k}{d \nu_k}}
d \nu_k >0 \,\,\, \Leftrightarrow \,\,\, \prod_{k=1}^\infty \left(
\sum_{i \in \mathbb{N}_k } \sqrt{p_{ik} q_{ik}}\right)>0,
\label{eqv}\end{equation}
\begin{equation} \label{sing}
  \mu_\xi \perp \lambda \,\,\, \Leftrightarrow \,\,\, \prod_{k=1}^\infty \int_{\Omega_k} \sqrt{\frac{ d \mu_k}{d \nu_k}}
d \nu_k =0   \,\,\, \Leftrightarrow \,\,\, \prod_{k=1}^\infty
\left( \sum_{i \in \mathbb{N}_k } \sqrt{p_{ik} q_{ik}}\right)=0.
\end{equation}

Of course, a singularly distributed random variable $\xi$ can also
be distributed discretely. For any point $x \in [0,1]$ the set
$f^{-1} (x)$ consists of at most two points from $\Omega$.
Therefore, the measure $\mu_\xi$ is an atomic measure if and only if
the measure $\mu$ is atomic.


If $\prod\limits_{k=1}^\infty \max\limits_i p_{ik}=0$, then $$ \mu
(\omega) = \prod_{k=1}^\infty p_{\omega_k k} \le \prod_{k=1}^\infty
\max\limits_i p_{ik} =0 \,\,\, \mbox{for any } \, \omega \in
\Omega,$$ and $\mu$ is continuous.

 If $\prod\limits_{k=1}^\infty
\max\limits_i p_{ik}>0$, then we consider the subset $A_+ =\{ \omega
: \, \mu (\omega)>0 \}$. The set $A_+$ contains the point
$\omega^*=(\omega^*_1, \omega^*_2,~...,~\omega^*_k,~...)$ such that
$p_{\omega^*_k k} = \max\limits_i p_{ik}$. It is easy to see that
for all $\omega \in A_+$ the condition $p_{\omega_k k} \not=
\max\limits_i p_{ik}$ holds only for a finite amount of values $k$.
This means that $A_+$ is a countable set and the event "$\omega \in
A_+$'' does not depend on any finite coordinates of $\omega$.
Therefore, by using Kolmogorov's ''$0$ and $1$'' theorem, we
conclude that $\mu (A_+) =0$ or $\mu (A_+)=1$. Since $\mu(A_+) \ge
\mu (\omega^*) >0$, we have $\mu(A_+) =1$, which proves the
discreteness of the measure $\mu$. \endproof
\textbf{Remark 3.} If
there exists a positive number $q^{+}$ such that $q_{ik}\geq
q^{+},\forall k\in \mathbf{N},$  $ \forall i\in \mathbf{N}_k,$ then
condition (\ref{eqv}) is equivalent to the convergence of the
following series:

\begin{equation}
\sum^{\infty }_{k=1}\{ \sum _{i\in \mathbf{N}_k}
(1-\frac{p_{ik}}{q_{ik}})^2 \} <\infty .  \label{Rem5}
\end{equation}

If $\lim\limits_{\overline{k \rightarrow \infty} } q_{ik} =0,$
then, generally speaking, conditions (\ref{eqv}) and (\ref
{Rem5}) are not equivalent. For example, let us consider the matrices $%
\widetilde{Q}\,$ and $\,\widetilde{P}\,$ as follows:
$\mathbf{N}_k=\left\{
0,1,2\right\} ,q_{1k}=\frac 1{2^k},\,\,\,\,q_{0k}=q_{2k}=\frac{1-q_{1k}}%
2,\,\,\,\,$ $p_{1k}=0,\,\,\,\,\,\,\,\,\,p_{0k}=p_{2k}=\frac 12.$
In this case condition (\ref{eqv}) holds, but (\ref{Rem5}) does
not hold.

\section{Metric-topological classification and fractal properties of the distributions of the random
variables with independent $\widetilde{Q}-$symbols}

For any probability distribution there exist sets which
essentially characterize the properties of the distribution. We
would like to stress the role of the following sets.

a) Topological support $S_\psi = \{ x: \, F(x+ \varepsilon) - F(x-
\varepsilon ) > \varepsilon, \,\, \forall \varepsilon>0\}$.
$S_\psi$ is the smallest closed support of the distribution of
$\psi$.

b) Essential support $N_\psi^\infty = \left\{ x: \,\,
\lim\limits_{\overline{\varepsilon \to 0}} \frac{F(x+ \varepsilon)
- F(x- \varepsilon )}{2\varepsilon} =+ \infty \right\}.$

 If the
topological support of a distribution is a fractal, then the
corresponding distribution is said to be externally fractal. The
probability distribution of a random variable $\psi$ is said to be
internally fractal if the essential support of the distribution is
a fractal set.

First of all we shall analyze the metric and topological properties
of the topological support of the random variable with independent
$\widetilde{Q}$-symbol.  In \cite{AKT, P} it was proven that
arbitrary singular continuous probability measures can be decomposed
into linear combinations of singular probability measures of S-, C-
and P-types (see the Introduction for the definitions).

We shall prove now that the above considered probability measures
$\mu_{\xi } $ are of the pure above mentioned metric-topological
types. Moreover we give necessary and sufficient conditions for a
probability measure to belong to each of these types.

\textbf{Theorem 3}. \textit{The distribution of the random
variable $\xi $ with independent $\widetilde{Q}-$symbols has pure
metric-topological type. Namely, the measure $\mu _\xi $ is one of
the  following three types: }

\textit{1) it is of the pure S-type if and only if the matrix
$\widetilde{P} $ contains only a finite number of columns
containing zero elements; }

\textit{2) it is of the pure C-type if and only if the matrix
$\widetilde{P} $ contains infinitely many columns having some
elements $p_{ik}=0$, and besides
\begin{equation}
\sum_{k=1}^\infty (\sum_{i:p_{ik}=0}q_{ik})=\infty~~~.  \label{Ctyp}
\end{equation}
}

\textit{3) it is of the pure P-type if and only if the matrix
$\widetilde{P} $ contains infinitely many columns having zero
elements and besides
\begin{equation}
\sum_{k=1}^\infty ( \sum_{i: p_{ik}=0} q_{ik}) <\infty ;  \label{Ptyp}
\end{equation}
}

\proof Let us consider the set $\Gamma \equiv \Gamma _{\widetilde{Q}%
(\mathbf{V)}}$ (see Sect. 3) with
$\mathbf{V}=\{\mathbf{V}_k\}_{k=1}^\infty $ defined by the
$\widetilde{P}-$matrix as follows: $\mathbf{V}_k=\{i\in
\mathbf{N}_k:~p_{ik}\not= 0\}$. It is easy to see that the
topological support of the measure ${\mu _\xi }$ coincides with a
set $\Gamma $ or its closure, i.e.,
\begin{equation}
S_{\xi }=\Gamma _{\widetilde{Q}(\mathbf{V)}}^{cl}.
\label{sigmaGamma}
\end{equation}
Therefore to examine the metric-topological structure of the set
$S_{\xi }$ we may apply the results of section 3. So, if the matrix
$\widetilde{P}$ contains only finite number of zero elements, then
$\mathbf{V}_k=\mathbf{N}_k,\ k>k_0$ for some $k_0>0$. In such a
case, $\Gamma $ is a union of at most  countable family of closed
intervals. Hence (\ref{sigmaGamma}) implies that the measure $\mu
_\xi $ is of the pure S-type.

In the opposite case where the matrix $\widetilde{P}$ contains an
infinite number of columns where some elements $p_{ik}=0$, then
obviously $\Gamma$ is a nowhere dense set (see Sect. 3). The
Lebesgue measure of the set $\Gamma $ by Lemma 1 is equal to
\begin{equation*}
\lambda (\Gamma )=\prod_{k=1}^\infty S_k(\mathbf{V}%
)=\prod_{k=1}^\infty (\sum_{i\in \mathbf{V}_k}q_{ik})=\prod_{k=1}^\infty
(1-\sum_{i:p_{ik}=0}q_{ik}).
\end{equation*}

 Then, by Lemma 2, either $\lambda(
\Gamma) =0$, provided that
condition (\ref{Ctyp}) fulfilled, or $\lambda( \Gamma ) >0$, if condition (%
\ref{Ptyp}) holds. Thus the measure $\mu _\xi $ either is of the
C-type, or it is of the P-type.

Since the conditions 1), 2) and 3) of this theorem are mutually
exclusive and one of them always holds, we conclude that the
distribution of the random variable $\xi $ with independent
$\widetilde{Q}$-symbols always has a pure metric-topological type.
\endproof
By using the latter theorems we can construct measures of
\textbf{8} kinds: pure point as well as pure singular continuous
of any S-, C-, or P-types, and  pure absolutely continuous
 but only of the S- and P-types.

We illustrate this statement by examples.

\textbf{Example 3}.

Let $\mathbf{N}_k=\{0,1,2\}$ and let the $\widetilde{Q}-$matrix be given by $%
q_{0k}=q_{1k}=q_{2k}=\frac 13,\ k\in N$.

S$_{pp}$: If $p_{0k}=\frac{1-p_{1k}}2, p_{1k}=1-\frac 1{2^k},p_{2k}= \frac{%
1-p_{1k}}2,$ then $\mu _\xi $ is a discrete measure of the pure
S-type. In this case $S _{\xi }= [ 0,1] $ and $N^{\infty} _{\xi }$
is a countable set which is dense on $[ 0,1] $.

S$_{sc}$: If $p_{0k}=\frac 14, p_{1k}=\frac 12, p_{2k}=\frac 14,$ then $%
\mu_\xi $ is a singular continuous measure of pure S-type. In this
case again $S _{\xi }= [ 0,1] $ but $N^{\infty} _{\xi }$ is now a
fractal set which is also dense on $[ 0,1] $.

S$_{ac}$: If $p_{0k}=p_{1k}=p_{2k}=\frac 13,$ then $\mu_\xi $
coincides with  Lebesgue measure on $[ 0,1] $.
 \vskip 0.3cm
\textbf{Example 4}.

Let again $\mathbf{N}_k=\{ 0,1,2\}$ and let the
$\widetilde{Q}-$matrix be given by $q_{0k}=q_{1k}=q_{2k}=\frac 13
, \ k\in N$. Then

C$_{pp}$: If $p_{0k}=1 - \frac{1}{2^k}, p_{1k}=0, p_{2k}=\frac
1{2^k}$, then $\mu _\xi $ is a pure point measure of the pure
C-type. In this case $S_{\xi }\equiv C_0 $ coincides with the
classical Cantor set $C_0$ and its essential support is a
countable set which is dense on $C_0 $.

C$_{sc}$: If $p_{0k}=\frac 1{2}, \,p_{1k}=0, p_{2k}=\frac1{2}$, then $%
\mu_\xi $ is the classical Cantor measure.


\vskip 0.3cm
 \textbf{Example 5}.

Let as above $\mathbf{N}_k=\{ 0,1,2\}$ and let the
$\widetilde{Q}-$matrix be given by
$q_{0k}=q_{2k}=\frac{1-q_{1k}}{2}, q_{1k}=\frac 1{2^k} , \ k\in N
$ Then

P$_{pp}$: If $p_{0k}=1-\frac 1{2^k}, p_{1k}=0, p_{2k}=\frac 1{2^k},$ then $%
\mu_\xi $ is a pure point measure of the pure P-type.

P$_{sc}$: If $p_{0k}=\frac 14, p_{1k}=0, p_{2k}=\frac 34,$ then
$\mu_\xi $ is a singular continuous measure of the pure P-type.

P$_{ac}$: If $p_{0k}=p_{2k}=\frac{1}2, p_{1k}=0,$ then $%
\mu_\xi $ is of the pure P-type measure which is absolutely
continuous w.r.t. Lebesgue measure.

We would like to stress that the essential support is more suitable
to describe the properties of distributions with complicated  local
structure. As we saw above, a discrete probability distribution may
be of $C$-, $P$- resp. $S$-type, and the topological support of
discrete distribution can be of any Hausdorff-Besicovitch dimension
$\alpha_0 \in [0,1]$. But the essential support of a discrete
distribution is always at most countable set.

The essential support is especially suitable for singular
distributions because of the following fact: a random variable
$\psi$ is singularly distributed iff $P_\psi (N_\psi^\infty)=1$.

For an absolutely continuous distribution the topological support
is always of positive Lebesgue measure. But the essential support
may be of very complicated local structure. In \cite{APT} we
constructed an example of an absolutely continuous distribution
function such that the essential support is an everywhere dense
superfractal set ($\alpha_0 (N_\xi^\infty) =1$). Therefore, the
condition $\alpha_0 (N_\xi^\infty) >0$ does not imply the
singularity of the distribution.

The following notion is very important for describing the fractal
properties of probability distributions. Let $A_\xi$ be the set of
all possible supports of the distribution of the r.v. $\xi$, i.e.,
$$ A_\xi = \{ E: \,\, E \in \mathcal{B}, \,\, P_\xi (E) =1 \}.$$
The number $\alpha_0 (\xi) = \inf\limits_{E \in A_\xi} \{ \alpha_0
(E) \}$ is said to be the Hausdorff-Besicovitch dimension of the
 distribution of the r.v. $\xi$.

It is obvious that $\alpha_0 (\xi) =0$ for any discrete
distribution; on the other hand,  $\alpha_0 (\xi) =1$ for any
absolutely continuous distribution.  $\alpha_0 (\xi)$ can be an
arbitrary number from $[0,1]$ for a singular continuous
distribution.

The determination of the Hausdorff-Besicovitch dimension of the
distribution of a random variable $\psi$ is essentially more
difficult problem  than the calculation of the fractal dimension
of the corresponding topological support, and in general setting
the problem of the determination of the Hausdorff-Besicovitch
dimension of the distribution of the random variable $\xi$ with
independent $\widetilde{Q}$-symbols is still open. The next (and
more complicated) stage in the fractal analysis of  a singularly
distributed random variable $\psi$  is the finding the support
$E_{\psi}
 $ of $\psi$ such that $\alpha_0(E_{\psi})=\alpha_0(\psi).$

{\bf Theorem 4.} {\it Let $p_{ik}= p_i$, $q_{ik}= q_i$, $ k \in N$,
$i \in N^0_{s-1}$, and let
$$
M[\widetilde{Q},(p_0, ..., p_{s-1})] = \left\{ x: \lim\limits_{k
\to \infty} \frac{N_i (x,k)}{k} = p_i,  i \in N^0_{s-1} \right\},
$$
where $N_i (x,k) $ is the number  of symbols ''$i$'' among the
first $k$ symbols of the $\widetilde{Q}$-representation of $x$.

Then $E_{\xi}=M[\widetilde{Q}, (p_0, ..., p_{s-1})]$ and
  $$ \alpha_0 (\xi)=
\frac{\sum\limits_{i=0}^{s-1} p_i \ln p_i}{
\sum\limits_{i=0}^{s-1} p_i \ln q_i}. \eqno  $$}

\proof For the simple $Q-$ representation the problem of the
determination  of the Hausdorff-Besicovitch dimension of the
distribution of the random variable $\xi$ was solved in \cite{Tor}.
In particular, we have $ \alpha_0 (\xi)=
\frac{\sum\limits_{i=0}^{s-1} p_i \ln p_i}{ \sum\limits_{i=0}^{s-1}
p_i \ln q_i}.$ From Theorem 1 it follows that the
Hausdorff-Besicovitch dimension of the  set $M [\widetilde{Q},
(p_0,...,p_{s-1})]$, which consists of the points whose
$\widetilde{Q}$-representation contains the digit ''$i$'' with the
asymptotic frequency $p_i$, is equal to
$\frac{\sum\limits_{i=0}^{s-1} p_i \ln p_i}{ \sum\limits_{i=0}^{s-1}
p_i \ln q_i}.$ So, the set $M [\widetilde{Q}, (p_0,...,p_{s-1})]$
can be considered as the ''dimensionally minimal'' support of the
distribution of the random variable with independent identically
distributed $Q$-symbols
\endproof

\textbf{Acknowledgement}

This work was partly supported by DFG 436 UKR 113/80 and DFG 436 UKR
113/67, INTAS 00-257, and SFB-611 projects. The authors would like
to express their gratitude to the referees for helpful remarks and
comments. The last three named authors gratefully acknowledge the
hospitality of the Institute of Applied Mathematics and of the IZKS
of the University of Bonn.

\end{document}